\documentclass[11pt]{article}
\usepackage{mathrsfs}
\usepackage{amsthm}
\usepackage{amssymb}
\usepackage{amsmath}
\usepackage{graphicx}
\usepackage{svg}
\usepackage{color}
\usepackage{amsfonts}
\usepackage{float}
\usepackage{cite}
\usepackage[text={140mm,210mm},left=35mm,vmarginratio=1:1]{geometry}
\newtheorem{theorem}{Theorem}[section]

\newtheorem{lemma}[theorem]{Lemma}
\newtheorem{proposition}[theorem]{Proposition}

\newtheorem{assumption}[theorem]{Assumption}
\numberwithin{equation}{section}
\normalsize

\begin{document}
\title{\textbf{Moderate Deviations of Hitting Times for Trajectories of Sums 
		of Independent and Identically Distributed Random Variables}}

\author{Yuheng He \thanks{\textbf{E-mail}: heyuheng@bjtu.edu.cn \textbf{Address}: School of Mathematics and Statistics, Beijing Jiaotong University, Beijing 100044, China.}
\\ Beijing Jiaotong University}

\date{}
\maketitle

\noindent {\bf Abstract:} In this paper we establish a moderate deviation principle of the hitting times for trajectories of sums 
of independent and identically distributed random variables. The main idea of proof is to convert the moderate deviations over a small time interval into the moderate deviations at a point, then utilize the moderate deviations for trajectories of sums
of independent random variables given by Hu in \cite{Hu2001} to get the moderate deviations at a point. By the upper bounds of large deviations for trajectories of sums 
of independent random variables given by Schuette in \cite{ Schuette1994}, we can prove that our convert doesn't influence the rate function of the moderate deviation.

\quad


\section{Introduction and main results}\label{section one}

Given a sequence of i.i.d. real valued random variables $\{X_i\}_{i=1}^{\infty}$,  satisfying $E(X_1)=\mu>0$ and ${\rm Var}(X_1)=\sigma^2$. The trajectories $\widetilde S_n$ of  $\{X_i\}_{i=1}^{\infty}$ studied in this paper are defined as follows
\[
\widetilde S_n(t)=\sum_{i=1}^{\lfloor nt \rfloor}X_i+(nt-{\lfloor nt \rfloor})X_{{\lfloor nt \rfloor}+1}
\]
for any $t\in[0,1]$, where $\lfloor a \rfloor$ represents the largest integer that does not exceed $a$. In this paper, we always assume that
\begin{assumption}\label{assumption 1}
	\[
\Lambda(a)={\rm log}E{\rm e}^{aX_1}<\infty  
    \]
	for all  $a\in\mathbb{R}$, $\Lambda(a)$ is continuous, and
 \[
\sup_{a\in\mathbb{R}}\{-\Lambda(a)\}<\infty.  
    \]

\end{assumption}

\begin{assumption}\label{assumption 2}
	There exists some $\theta\in(0,1],~v>1$, and $b>0$ such that 
 \[
 E{\rm exp}(\theta|X_1|^v)\leq{\rm e}^b.
 \]

\end{assumption}

Using the Law of Large Numbers and converging together lemma given in \cite{Durrett2019}, it is easy to see that
\[
\lim_{n\rightarrow+\infty}\frac{\widetilde S_n(t)}{n}=x(t)
\]
in probability for any $t\in[0,1]$, where $x(t) = \mu t$. Since $\{x(t)\}_{0\leq t\leq 1}$ is increasing with $t$, there is one and only one $\tau_r>0$ satisfying equation $x(\tau_r)=r$ for any $r\in (0, \mu)$. It's easy to verify that  
\begin{equation}\label{equ tau}
\tau_r = \frac{r}{\mu}.
\end{equation}

Similarly we define the hitting times $\tau_r^n$ for trajectories of sums, for any $r\in (0, \mu)$, let
\[
\tau_r^n=\inf\left\{t:~\widetilde S_n(t) \geq nr\right\}.
\]

Using the continuous mapping theorem in \cite{Billingsley1969} of Billingsley and central limit theorem for trajectories in \cite{Donsker1951} of Donsker together with a similar method as Section 11.4 in \cite{Ethier1986} of Ethier and Kurtz, we can derive the following propositions.

\begin{proposition}

1)
\[
\lim_{n\rightarrow+\infty}\tau_r^n=\tau_r
\]
in probability for any $r\in (0, \mu)$.

2)
\[
\sqrt{n}\left(\tau_r^n-\tau_r\right) \Rightarrow_n \frac{\sigma}{\mu}B_{\tau_r}
\]
for any $r\in (0, \mu)$, where $B_t\sim\mathcal{N}(0, t)$.
\end{proposition}

Since we have already established the Law of Large Numbers and the Central Limit Theorem of $\tau_r^n$, we further consider deriving the moderate deviation principle of $\tau_r^n$ as the main theorem of this paper.

\begin{theorem}\label{theorem main moderate deviation}
Let $\{a_n\}_{n=1}^{\infty}$ be a positive sequence satisfying that
\[
\lim_{n\rightarrow+\infty}\frac{a_n}{n}=\lim_{n\rightarrow+\infty}\frac{\sqrt{n}}{a_n}=0,
\]
then
\begin{equation}\label{equ 1.1}
\lim_{n\rightarrow+\infty}\frac{n}{a_n^2}\log P\left(\frac{n}{a_n}(\tau_r^n-\tau_r)>t\right)=-\frac{\mu^3t^2}{2\sigma^2r}
\end{equation}
and
\begin{equation}\label{equ 1.2}
\lim_{n\rightarrow+\infty}\frac{n}{a_n^2}\log P\left(\frac{n}{a_n}(\tau_r^n-\tau_r)<-t\right)=-\frac{\mu^3t^2}{2\sigma^2r}
\end{equation}
for any $r\in (0, \mu)$ and $t>0$.
\end{theorem}

We will show the proof of Theorem \ref{theorem main moderate deviation} in detail in Section \ref{section proof}. The main idea of the proof is to convert the events that occur during a small time interval into events that occur within closed or open sets, so that we can apply the results of moderate deviations for trajectories of sums
of independent random variables of Hu given in \cite{Hu2001}. Utilizing the upper bound of the large deviation for trajectories of sums of independent
random variables of Schuette given in \cite{ Schuette1994}, we can show the convert does not affect the rate function of the main event.

\section{Preliminary results}\label{section Preliminaries}
In this section we will prepare some key results which will play a significant role in the proof of the main theorem. We review moderate deviations for trajectories of sums
of independent random variables and provide an upper bound of the large deviation of $\tau_r^n$.

\subsection{A review of moderate deviations for trajectories of sums
	of independent random variables}
For any $t\geq 0$, let $\theta_t^n=\frac{\widetilde S_n(t)-nx(t)}{a_n}$, then for any given $T > 0$,  $\theta^{n,T}:=\{\theta_t^n\}_{0\leq t\leq T}$ is a random element in $\mathcal{C}\left([0, T], \mathbb{R}\right)$, where $\mathcal{C}\left([0, T], \mathbb{R}\right)$ represents the set of all continuous functions from $[0, T]$ to $\mathbb{R}$. Trajectories of sums
of independent random variables in this paper are those talked in \cite{Hu2001} under our assumptions in the beginning. Therefore, by Theorem 1.2 of \cite{Hu2001}, we obtain the following proposition.

\begin{proposition}\label{Proposition 2.1}(Hu, 2001, \cite{Hu2001}, Theorem 1.2)Under Assumptions \ref{assumption 1} and \ref{assumption 2},
for every $f\in \mathcal{C}\left([0, T], \mathbb{R}\right)$, let
\[
I_T(f)=
\begin{cases}
& \frac{1}{2\sigma^2}\int_0^T\left(f^\prime(u)\right)^2du \text{\quad if~}f\text{~is absolutely continuous and~}f(0)=0,\\
& +\infty \text{\quad otherwise},
\end{cases}
\]
then
\[
\limsup_{n\rightarrow+\infty}\frac{n}{a_n^2}\log P\left(\theta^{n,T}\in F\right)\leq -\inf_{f\in F}I_T(f)
\]
for every closed set $F\subseteq \mathcal{C}\left([0, T], \mathbb{R}\right)$ and
\[
\liminf_{n\rightarrow+\infty}\frac{n}{a_n^2}\log P\left(\theta^{n,T}\in G\right)\geq -\inf_{f\in G}I_T(f)
\]
for every open set $G\subseteq \mathcal{C}\left([0, T], \mathbb{R}\right)$.
\end{proposition}

The proof of Theorem \ref{theorem main moderate deviation} relies heavily on the following property of $I_T$.

\begin{proposition}\label{proposition 2.2 lower bound of IT}

\[
\inf\left\{I_T(f):~f\in \mathcal{C}\left([0, T], \mathbb{R}\right)\text{~and~}f(T)=a\right\}=\frac{a^2}{2\sigma^2T}
\]
for any $a\in \mathbb{R}$.
\end{proposition}

\proof[Proof of Proposition \ref{proposition 2.2 lower bound of IT}]

For any absolutely continuous $f$ with $f(0)=0$ and $f(T)=a$, by Cauchy-Schwarz's inequality,
\begin{align*}
a^2=f^2(T)=\left(\int_0^Tf^\prime(u)du\right)^2 
\leq \int_0^T\left(f^\prime(u)\right)^2du\int_0^T1^2du=2\sigma^2TI_T(f)
\end{align*}
therefore
\[
I_T(f)\geq \frac{a^2}{2\sigma^2T}.
\]
However, if we take $f^{*}(t)=\frac{a}{T}t$, then $f^{*}(0)=0, f^{*}(T)=a$ and
\[
I(f^{*})=\frac{1}{2\sigma^2}\int_0^T \left(\frac{a}{T}\right)^2du=\frac{a^2}{2\sigma^2T}
\]
 hence the proof is complete.

\qed

\subsection{An upper bound of the large deviation of the hitting time}

We give an upper bound of the large deviation of $\tau_r^n$ in this subsection.

\begin{lemma}\label{lemma 2.3 LDPupperboundHittingTime}
\[
\limsup_{n\rightarrow+\infty}\frac{1}{n}\log P\left(|\tau_r^n-\tau_r|>\epsilon
\right)<0
\]
for any $r\in (0, \mu)$ and $\epsilon>0$.
\end{lemma}
According to Lemma \ref{lemma 2.3 LDPupperboundHittingTime}, $P\left(|\tau_r^n-\tau_r|\leq\epsilon\right)$ converges to $1$ with an exponential speed as $n\rightarrow+\infty$, which 
shows our convert does not affect the rate function of moderate deviation in the proof of Theorem \ref{theorem main moderate deviation}.

In the proof of Lemma \ref{lemma 2.3 LDPupperboundHittingTime} we utilize upper bounds of large deviations for trajectories of sums of independent random variables given in \cite{ Schuette1994}. Under our assumptions, trajectories of sums of independent and identically distributed random variables in this paper are the special cases of those talked in \cite{ Schuette1994}. Therefore, we have the following the proposition by Theorem 3.1 of \cite{ Schuette1994}.

\begin{proposition}\label{Proposition 2.4 UpperBoundofLargeDeviation}(Schuette, 1994, \cite{ Schuette1994}, Theorem 3.1)
Under Assumptions \ref{assumption 1} and \ref{assumption 2}, for any $f\in \mathcal{C}\left([0, T], \mathbb{R}\right)$, let
\[
J_T(f)=
\begin{cases}
& \int_0^T\Lambda^{*}(f^\prime(u))du \text{\quad if~}f\text{~is absolutely continuous and~}f(0)=0,\\
& +\infty \text{\quad otherwise},
\end{cases}
\]
where
\[
\Lambda^{*}(x)=\sup_{a\in \mathbb{R}}\left\{a x-{\rm log}E{\rm e}^{a X_1}\right\},
\]
then $J_T$ is a good rate function and
\[
\limsup_{n\rightarrow+\infty}\frac{1}{n}\log P\left(\left\{\frac{\widetilde S_n(t)}{n}\right\}_{0\leq t\leq T}\in F\right)\leq -\inf_{f\in F}J_T(f)
\]
for every closed $F\subseteq \mathcal{C}\left([0, T], \mathbb{R}\right)$.
\end{proposition}

Notice that $\Lambda^{*}(x)\geq 0$ for any $x\in\mathbb{R}$, since let $\mathcal{R}_{x}(a):=ax-{\rm log}E{\rm e}^{a X_1}$, and when $a=0$, $\mathcal{R}_{x}(0)=0$, therefore, $J_T(f)\geq 0$. Furthermore, $J_T(f)=0$ if and only if $f(u)=x(u)$ for any $0\leq u\leq T$. Since $\partial_{aa}\mathcal{R}_{f^\prime(u)}(a)\leq 0$ by Cauchy-Schwarz's inequality, $\mathcal{R}_{f^\prime(u)}(a)$ is a concave function which has one and only one maximum.  Hence $\max_{a} \mathcal{R}_{f^\prime(u)}(a)=\mathcal{R}_{f^\prime(u)}(0)$ if and only if $\partial_a\mathcal{R}_{f^\prime(u)}(0)=0$, thus $f^\prime(u)=\mu$ and $f(u)=\mu u$ since $f(0)=0$.

Now we prove Lemma \ref{lemma 2.3 LDPupperboundHittingTime}.

\proof[Proof of Lemma \ref{lemma 2.3 LDPupperboundHittingTime}]

Let $0<\epsilon<\tau_r$, then
\[
\left\{\tau_r^n-\tau_r<-\epsilon\right\}\subseteq \left\{\left\{\frac{\widetilde S_n(t)}{n}\right\}_{0\leq t\leq \tau_r-\epsilon}\in A_\epsilon^-\right\},
\]
where
\[
A_\epsilon^-=\left\{f\in \mathcal{C}\left([0, \tau_r-\epsilon], \mathbb{R}\right):~\sup_{0\leq t\leq \tau_r-\epsilon}f(t)\geq r\right\}.
\]
It is easy to verify that $A_\epsilon^-$ is closed in $\mathcal{C}\left([0, \tau_r-\epsilon], \mathbb{R}\right)$ by using the properties of uniformly continuous functions. Notice that $\{x(t)\}_{0\leq t\leq \tau_r-\epsilon}\not\in A_\epsilon^-$ since $\sup_{0\leq t\leq \tau_r-\epsilon}x(t) = \mu (\tau_r-\epsilon)=r-\mu\epsilon< r$.  $J_{\tau_r-\epsilon}(f)>0$ for $\{f(t)\}_{0\leq t\leq \tau_r-\epsilon}\neq \{x(t)\}_{0\leq t\leq \tau_r-\epsilon}$ from the analysis above, and $J_{\tau_r-\epsilon}$ is a good rate function by Proposition \ref{Proposition 2.4 UpperBoundofLargeDeviation}, we have
\[
\inf_{f\in A_\epsilon^-}J_{\tau_r-\epsilon}(f)>0
\]
and hence
\begin{equation}\label{equ 2.1}
\limsup_{n\rightarrow+\infty}\frac{1}{n}\log P\left(\tau_r^n-\tau_r<-\epsilon\right)\leq -\inf_{f\in A_\epsilon^-}J_{\tau_r-\epsilon}(f)<0.
\end{equation}
 On the other hand,

\[
\left\{\tau_r^n-\tau_r>\epsilon\right\}\subseteq \left\{\left\{\frac{\widetilde S_n(t)}{n}\right\}_{0\leq t\leq \tau_r+\epsilon}\in A_\epsilon^+\right\},
\]
 
where
\[
A_\epsilon^+=\left\{f\in \mathcal{C}\left([0, \tau_r+\epsilon], \mathbb{R}\right):~\sup_{0\leq t\leq \tau_r+\epsilon}f(t)\leq r\right\}.
\]
Similarly, it is easy to check that $A_\epsilon^+$ is also closed and $\{x(t)\}_{0\leq t\leq \tau_r+\epsilon}\not\in A_\epsilon^+$. Therefore,
\[
\inf_{f\in A_\epsilon^+}J_{\tau_r+\epsilon}(f)>0
\]
and hence
\begin{equation}\label{equ 2.2}
\limsup_{n\rightarrow+\infty}\frac{1}{n}\log P\left(\tau_r^n-\tau_r>\epsilon\right)\leq -\inf_{f\in A_\epsilon^+}J_{\tau_r+\epsilon}(f)
<0.
\end{equation}
 By Equation \eqref{equ 2.1} and \eqref{equ 2.2}, the proof of Lemma \ref{lemma 2.3 LDPupperboundHittingTime} is complete.

\qed

\section{The proof of our main theorem}\label{section proof}

In this section, we show the proof of Theorem \ref{theorem main moderate deviation}. 
Since Equation \eqref{equ 1.1} and \eqref{equ 1.2} can be proved analogously, we only show the proof of Equation \eqref{equ 1.1} in detail.
To prove Equation \eqref{equ 1.1}, we consider proving the following two equations
\begin{equation}\label{equ 3.1}
\limsup_{n\rightarrow+\infty}\frac{n}{a_n^2}\log P\left(\frac{n}{a_n}(\tau_r^n-\tau_r)>t\right)\leq-\frac{\mu^3t^2}{2\sigma^2r}
\end{equation}
and
\begin{equation}\label{equ 3.2}
\liminf_{n\rightarrow+\infty}\frac{n}{a_n^2}\log P\left(\frac{n}{a_n}(\tau_r^n-\tau_r)>t\right)\geq-\frac{\mu^3t^2}{2\sigma^2r}.
\end{equation}
First we prove Equation \eqref{equ 3.1}.

\proof[Proof of Equation \eqref{equ 3.1}]

By the definition of $\tau^n_r$, we have $\widetilde S_n(\tau^n_r)=nr$ and notice that $x(t)=\mu t$. Hence,
\begin{equation}\label{equ 3.3}
\left\{\frac{n}{a_n}(\tau_r^n-\tau_r)>t\right\}
= \left\{\frac{\widetilde S_n(\tau^n_r)-nx(\tau^n_r)}{a_n}<-\mu t\right\}.
\end{equation}

For any $0<\delta<\tau_r$,
\[
\left\{\frac{\widetilde S_n(\tau^n_r)-nx(\tau^n_r)}{a_n}<-\mu t, |\tau_r^n-\tau_r|\leq \delta\right\}
\subseteq \left\{\theta^{n,\tau_r+\delta}\in F_{ \delta}\right\},
\]
where $\theta^{n, \tau_r+\delta}=\{\frac{\widetilde S_n(t)-nx(t)}{a_n}\}_{0\leq t\leq \tau_r+\delta}$ as defined in Section \ref{section Preliminaries} and
\[
F_{\delta}=\left\{f\in \mathcal{C}\left([0, \tau_r+\delta], \mathbb{R}\right):~\inf_{\tau_r-\delta\leq s\leq \tau_r+\delta}f(s)\leq -\mu t\right\}.
\]
It's easy to check that $F_{\delta}$ is closed in $\mathcal{C}\left([0, \tau_r+\delta], \mathbb{R}\right)$, then by Equation \eqref{equ 3.3}, Lemma \ref{lemma 2.3 LDPupperboundHittingTime} together with Proposition \ref{Proposition 2.1},
\begin{equation}\label{equ 3.4}
\limsup_{n\rightarrow+\infty}\frac{n}{a_n^2}\log P\left(\frac{n}{a_n}(\tau_r^n-\tau_r)>t\right)\leq -\inf_{f\in F_{\delta}}I_{\tau_r+\delta}(f).
\end{equation}
For any given $f\in F_{\delta}$, there exists $s(\delta)\in [\tau_r-\delta, \tau_r+\delta]$ such that $f(s(\delta))\leq -\mu t<0$. According to the definition of $I_T$ and Proposition \ref{proposition 2.2 lower bound of IT},
\begin{align*}
I_{\tau_r+\delta}(f)\geq I_{s(\delta)}\left(\{f(u)\}_{0\leq u\leq s(\delta)}\right)
\geq \frac{f^2(s(\delta))}{2\sigma^2s(\delta)}\geq \frac{\mu^2t^2}{2\sigma^2s(\delta)}.
\end{align*}
Hence
\[
\limsup_{n\rightarrow+\infty}\frac{n}{a_n^2}\log P\left(\frac{n}{a_n}(\tau_r^n-\tau_r)>t\right)\leq
-\frac{\mu^2t^2}{2\sigma^2s(\delta)}
\]
follows from Equation \eqref{equ 3.4}. Since $\delta$ is arbitrary, let it converge to $0$, then
\begin{equation}\label{equ 3.5}
\limsup_{n\rightarrow+\infty}\frac{n}{a_n^2}\log P\left(\frac{n}{a_n}(\tau_r^n-\tau_r)>t\right)\leq
-\frac{\mu^2t^2}{2\sigma^2\tau_r}=-\frac{\mu^3t^2}{2\sigma^2r}.
\end{equation}

The proof of Equation \eqref{equ 3.1} is complete.

\qed

Now we prove Equation \eqref{equ 3.2}.

\proof[Proof of Equation \eqref{equ 3.2}]

For any $0<\delta<\tau_r$,
\[
\left\{\frac{\widetilde S_n(\tau^n_r)-nx(\tau^n_r)}{a_n}<-\mu t, |\tau_r^n-\tau_r|\leq \delta\right\}
\supseteq \left\{\theta^{n,\tau_r+\delta}\in G_{\delta}, |\tau_r^n-\tau_r|\leq \delta\right\},
\]
where
\[
G_{\delta}=\left\{f\in \mathcal{C}\left([0, \tau_r+\delta], \mathbb{R}\right):~\sup_{\tau_r-\delta\leq s\leq \tau_r+\delta}f(s)<-\mu t\right\}.
\]
It's easy to verify that $G_{\delta}$ is open in $\mathcal{C}\left([0, \tau_r+\delta], \mathbb{R}\right)$, then by Equation \eqref{equ 3.3}, Lemma \ref{lemma 2.3 LDPupperboundHittingTime} together with Proposition \ref{Proposition 2.1},
\begin{equation}\label{equ 3.10}
\liminf_{n\rightarrow+\infty}\frac{n}{a_n^2}\log P\left(\frac{n}{a_n}(\tau_r^n-\tau_r)>t\right)\geq -\inf_{f\in G_{\delta}}I_{\tau_r+\delta}(f).
\end{equation}
For any $\epsilon>0$, let 
\[
f_{\epsilon, \delta}(s)=-\frac{\mu t(1+\epsilon)}{\tau_r-\delta}s,
\]
then
\[
\sup_{\tau_r-\delta\leq s\leq \tau_r+\delta}f_{\epsilon, \delta}(s)=-\mu t(1+\epsilon)<-\mu t,
\]
hence $f_{\epsilon, \delta}\in G_{\delta}$, and we have $I_{\tau_r+\delta}(f_{\epsilon, \delta}) = \frac{\mu^2t^2(\tau_r+\delta)(1+\epsilon)^2}{2\sigma^2(\tau_r-\delta)^2}$.
Then, by Equation \eqref{equ 3.10},
\[
\liminf_{n\rightarrow+\infty}\frac{n}{a_n^2}\log P\left(\frac{n}{a_n}(\tau_r^n-\tau_r)>t\right)\geq
-\frac{\mu^2t^2(\tau_r+\delta)(1+\epsilon)^2}{2\sigma^2(\tau_r-\delta)^2}.
\]
Since $\epsilon$ and $\delta$ are arbitrary, let them converge to $0$, then
\[
\liminf_{n\rightarrow+\infty}\frac{n}{a_n^2}\log P\left(\frac{n}{a_n}(\tau_r^n-\tau_r)>t\right)\geq 
-\frac{\mu^2t^2}{2\sigma^2\tau_r}=
-\frac{\mu^3t^2}{2\sigma^2r}.
\]
Hence Equation \eqref{equ 3.2} is proved.

\qed

\section{Applications}\label{section applications}

In this section we simulate Theorem \ref{theorem main moderate deviation} by Python in two examples. Throughout this section we take that $a_n=n^{0.9}$ which satisfies $\lim_{n\rightarrow+\infty}\frac{a_n}{n}=\lim_{n\rightarrow+\infty}\frac{\sqrt{n}}{a_n}=0$.

\textbf{Example 1} Let $X_i$ follow an exponential distribution with rate $\lambda=1.0$, then by Theorem \ref{theorem main moderate deviation} and Equation \eqref{equ tau} we obtain
\[
\lim_{n\rightarrow+\infty}\frac{n}{a_n^2}\log P\left(\frac{n}{a_n}(\tau_r^n-\tau_r)>t\right)=-\frac{t^2}{2r}
\]
for any $r\in(0,1)$ and $t>0$, where $\tau_r=r$.

\begin{figure}[H]
\centering
\includegraphics[scale=0.25]{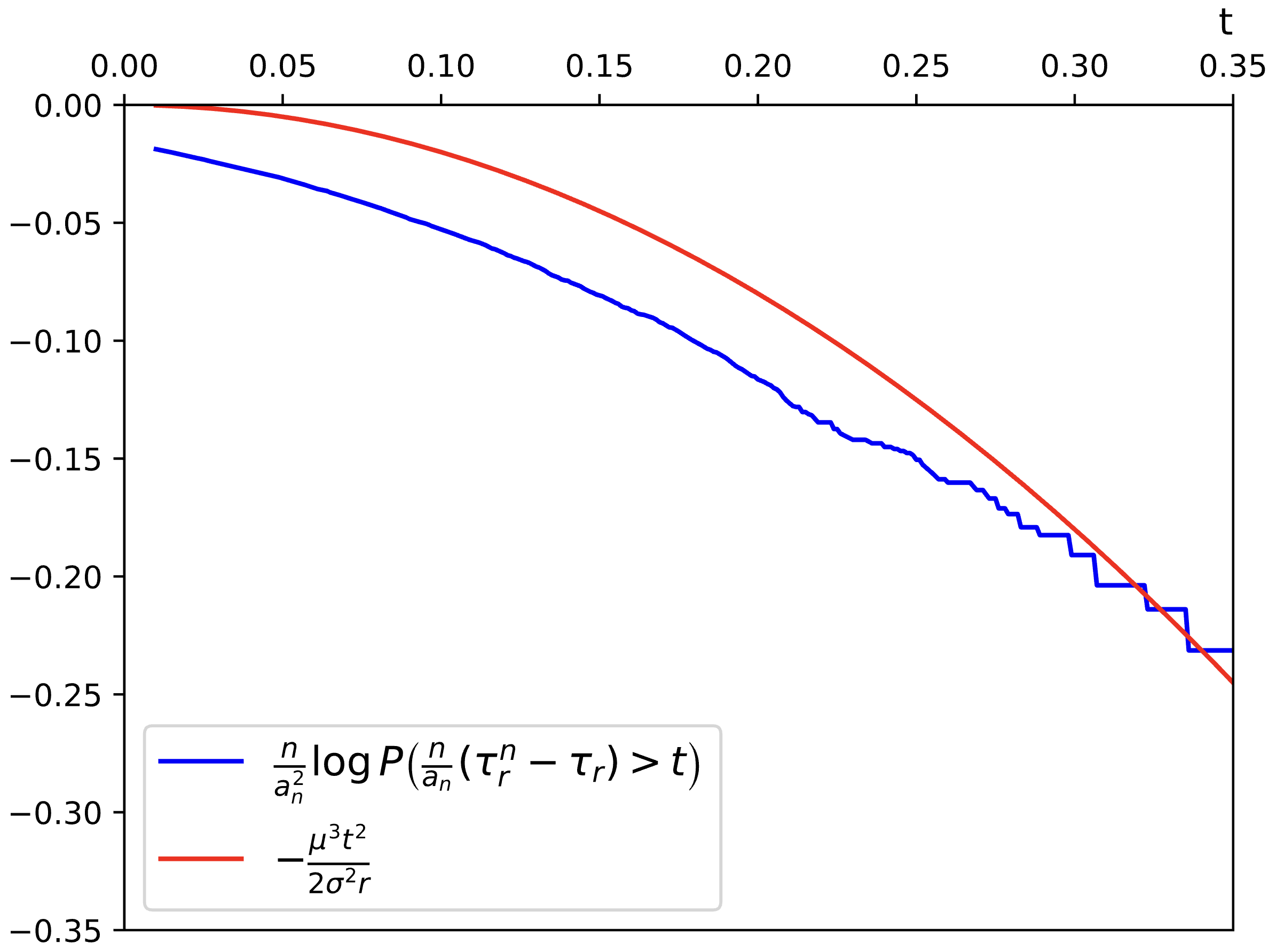}
\caption{$X_i\sim {\rm exponential}(1.0)$, $r=0.25$}\label{figure 4.1}
\end{figure}

We set $n=100$, $r=0.25$, and then simulate the mechanism of $\tau_{0.25}^{100}$. By a for-loop, we generate $10000$ independent copies of $\tau_{0.25}^{100}$. The blue polyline in Figure \ref{figure 4.1} shows results of
\[
-\frac{100}{(100^{0.9})^2}\log\left(\frac{\sum_{i=1}^{10000}1_{\{\frac{100}{100^{0.9}}(\tau_{0.25}^{100,i}-\tau_{0.25})>t\}}}{10000}\right)
\]
where $\tau_r^{100,i}$ is the $i^{\rm th}$ copy of $\tau_r^{100}$ and the red parabola shows our rate function $-2t^2$.

\textbf{Example 2} Let $X_i$ follow a Poisson distribution with $\lambda=1.0$, then by Theorem \ref{theorem main moderate deviation} and Equation \eqref{equ tau},
\[
\lim_{n\rightarrow+\infty}\frac{n}{a_n^2}\log P\left(\frac{n}{a_n}(\tau_r^n-\tau_r)>t\right)=-\frac{t^2}{2r}
\]
for any $r\in(0,1)$ and $t>0$, where $\tau_r=r$.

\begin{figure}[H]
\centering
\includegraphics[scale=0.25]{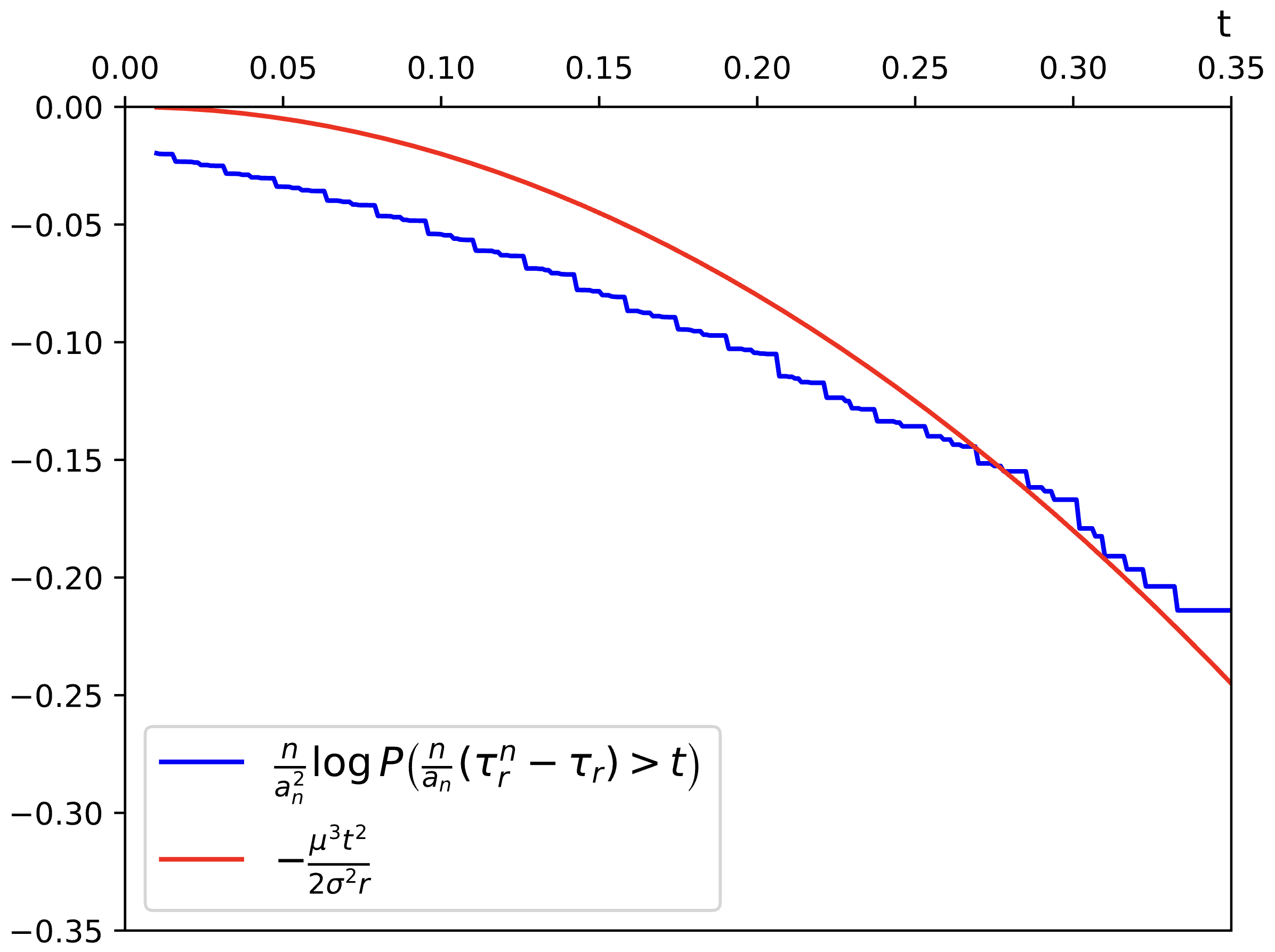}
\caption{$X_i\sim {\rm Poisson}(1.0)$, $r=0.25$}\label{figure 4.2}
\end{figure}

We also set $n=100$, $r=0.25$, and then simulate the mechanism of $\tau_{0.25}^{100}$. By a for-loop, we generate $10000$ independent copies of $\tau_{0.25}^{100}$. The blue polyline in Figure \ref{figure 4.2} shows results of
\[
-\frac{100}{(100^{0.9})^2}\log\left(\frac{\sum_{i=1}^{10000}1_{\{\frac{100}{100^{0.9}}(\tau_{0.25}^{100,i}-\tau_{0.25})>t\}}}{10000}\right)
\]
where $\tau_r^{100,i}$ is the $i^{\rm th}$ copy of $\tau_r^{100}$ and the red parabola shows our rate function $-2t^2$.

\quad

\textbf{Acknowledgments.} The author would like to thank the mentor and reviewers. Their comments are great help for the improvement of this paper. The authors are grateful to the financial support from the Fundamental Research Funds for the Central Universities with grant number 2022JBMC039.

{}
\end{document}